\documentclass[a4paper,11pt]{amsart}

\usepackage[english]{babel}
\usepackage{amsmath, enumerate, amsfonts, amssymb, amsthm}
\usepackage[dvips]{graphicx}

\newtheorem{defin}{\bf Def\mbox{}inition}
\newtheorem{theo}[defin]{\bf Theorem}
\newtheorem{prop}[defin]{\bf Proposition}

\newtheorem{cor}[defin]{\bf Corollary}

\newtheorem*{rem*}{\bf Remark}
\newtheorem*{nota*}{\bf Notation}

\newcommand{\dps}{\displaystyle}

\newcommand{\gauche}{\begin{flushleft}\end{flushleft}}

\newcommand{\C}{\mathbb{C}}
\newcommand{\N}{\mathbb{N}}

\newcommand{\Q}{\mathbb{Q}}
\renewcommand{\k}{\mathbf{k}}
\newcommand{\K}{\mathbf{K}}
\renewcommand{\L}{\mathbf{L}}

\newcommand{\CC}{\mathcal{C}}

\newcommand{\QQ}{\mathcal{Q}}
\newcommand{\PP}{\mathcal{P}}
\newcommand{\II}{\mathcal{I}}

\newcommand{\spec}{\text{Spec}}
\newcommand{\specm}{\text{Specm}}

\newcommand{\An} {\mathbf{A}_n}

\renewcommand{\O}{\mathcal{O}}

\newcommand{\D}{\mathcal{D}}

\newcommand{\dxi}{\partial _{x_i}}

\newcommand{\B}{\mathcal{B}}

\title[Generic Bernstein-Sato polynomial]{Generic
  Bernstein-Sato polynomial on an irreducible affine scheme}

\author{Rouchdi BAHLOUL}
\address{D\'epartement de Math\'ematiques\\ U.M.R. 6093\\Universit\'e
  d'Angers\\ 2 Bd Lavoisier, 49045 Angers cedex 01, France}
\email{rouchdi.bahloul@univ-angers.fr}

\begin{document}

\begin{abstract}
Given $p$ polynomials with coefficients in a commutative unitary
integral ring $\CC$ containing $\Q$, we define the notion of a generic
Bernstein-Sato polynomial on an irreducible affine scheme $V \subset
\spec(\CC)$. We prove the existence of such a non zero rational
polynomial which covers and generalizes previous existing results by
H. Biosca. When $\CC$ is the ring of an algebraic or analytic space,
we deduce a stratification of the space of the parameters such that on
each stratum, there is a non zero rational polynomial which is a
Bernstein-Sato polynomial for any point of the stratum. This
generalizes a result of A. Leykin obtained in the case $p=1$.
\end{abstract}

\maketitle

\section*{Introduction and Main Results}

Fix $n\ge 1$ and $p\ge 1$ two integers and
$v\in \N^p$. Let $x=(x_1,\ldots,x_n)$ and $s=(s_1,\ldots,s_p)$ be two
systems of variables. Let $\k$ be a field of 
characteristic \footnote{all the fields considered in this paper are
  of characteristic $0$} $0$. Let $\An(\k)$ be the ring of
differential operators with coefficients in $\k[x]=\k[x_1,\ldots,x_n]$
and $\D$ (resp. $\O$) be the sheaf of rings of differential operators
(resp. analytic functions) on $\C^n$ for which we denote by
$\D_{x_0}$ (resp. $\O_{x_0}$) the fiber in $x_0$.

Let $f=(f_1,\ldots, f_p)$ be in $\k[x]^p$ (resp. $\O_{x_0}^p$) and
consider the following functional identity:
\[b(s) f^s \in \An(\k)[s] \cdot f^{s+v},\]
(resp. $\D_{x_0}[s]$ instead of $\An(\k)[s]$) where
$f^{s+v}=f_1^{s_1+v_1} \cdots f_p^{s_p+v_p}$. This identity takes
place in the free module generated by $f^s$ over $\k[x,\frac{1}{f_1
  \cdots f_p},s]$ (resp. $\O_{x_0}[\frac{1}{f_1 \cdots f_p},s]$).

The set of such $b(s)$ is an ideal of $\k[s]$
(resp. $\C[s]$). This ideal is called the (global) Bernstein-Sato
ideal of $f$ (resp. local Bernstein-Sato ideal in $x_0$) and we denote
it by $\B^v(f)$ (resp. $\B^v_{x_0}(f)$). When $p=1$, this ideal is
principal and its monic generator is called the Bernstein polynomial
associated with $f$. Historically, I.N. Bernstein \cite{bernstein}
introduced the (global) Bernstein polynomial and proved its existence
(i.e. the fact that it is not zero). J.E. Bj\"ork \cite{bjorkPP} has
given the proof in the analytic case. Let us cite also M. Kashiwara
\cite{kashiwara} who proved, moreover, the rationality of the roots of
the local Bernstein polynomial. For $p\ge 2$, the algebraic case can
be easily treated in the same way as for $p=1$. For the analytic case,
the proof of the non nullity of $\B^v_{x_0}(f)$ is due to C. Sabbah
(\cite{sabbah1} and \cite{sabbah2}). Let us also cite A. Gyoja
\cite{gyoja} who proved that $\B^v_{x_0}(f)$ contains a non zero
\emph{rational} polynomial. The absolute Bernstein-Sato polynomial
naturally leads to the notion of a generic Bernstein-Sato polynomial
which we shall explain in what follows.\\

Let $\CC$ be a unitary commutative integral ring with the following
condition:\\
For any prime ideal $\PP \subset \CC$ and for any
$n\in \N\smallsetminus \{0\}$, we have:
\[n \in \PP \Rightarrow 1 \in \PP.\]
This condition is equivalent to the fact that for any $\PP\subset
\CC$, the fraction field of $\CC/\PP$ is of characteristic $0$. Note
that this condition is satisfied if and only if there exists an
injective ring morphism $\Q \hookrightarrow \CC$.

We shall see $\CC$ as the ring of coefficients or
parameters. Indeed, let $f=(f_1,\ldots, f_p)$ in $\CC[x]^p=
\CC[x_1,\ldots,x_n]^p$.

Let us denote by $\An(\CC)$ the ring of differential operators with
coefficients in $\CC[x]$, that is the $\CC$-algebra generated by $x_i$
and $\dxi$ ($i=1,\ldots,n$) where the only non trivial commutation
relations are $[\dxi, x_i]=1$ for $i=1,\ldots,n$ (hence $\CC$ is in
the center of $\An(\CC)$).

We denote by $\spec(\CC)$ (resp. $\specm(\CC)$) the set of prime
(resp. maximal) ideals of $\CC$ which is the spectrum of $\CC$
(resp. the maximal spectrum). For an ideal $\II\subset \CC$, we
denote by $V(\II)=\{\PP \in \spec(\CC)\, ;\, \PP \supset \II\}$ the
affine scheme defined by $\II$ and $V_m(\II)=V(\II) \cap
\specm(\CC)$. Remark that we shall only work with the closed subsets
of $\spec(\CC)$ and forget the sheaf structure of a scheme.

We are going to introduce the notion of a generic Bernstein-Sato
polynomial of $f$ on an irreducible affine scheme $V=V(\QQ) \subset
\spec(\CC)$ (that is when $\QQ$ is prime).\\
So let $\QQ$ be a prime ideal of $\CC$ and suppose that none of the
$f_j$'s is in $\QQ[x]$.\\
The main result of this article is the following.

\begin{theo}\label{theo:BSgen}
There exists $h \in \CC \smallsetminus \QQ$ and $b(s) \in
\Q[s_1,\ldots,s_p] \smallsetminus 0$ such that
\[h \, b(s)f^s \in \An(\CC)[s]f^{s+v} + \big( \QQ[x, \frac{1}{f_1
  \ldots f_p},s] \big) f^s.\]
\end{theo}
Such a $b(s)$ is called a (rational) generic Bernstein-Sato polynomial
of $f$ on $V=V(\QQ)$ (see the notation and the remark below).

In the case where $p=1$, the generic and relative (not introduced
here) Bernstein polynomial has been studied by F. Geandier in
\cite{geandier} and by J. Brian\c{c}on, F. Geandier and P. Maisonobe
in \cite{br-ge-m} in an analytic context (where $f$ is an analytic
function of $x$). In \cite{biosca} (see also
\cite{biosca2}), H. Biosca studied these notions with $p\ge 1$ in the
analytic and the algebraic context (that which we are concerned with)
and proved that when
\begin{itemize}
\item $\CC=\C[a_1,\ldots,a_m]$ or
\item $\CC=\C\{a_1,\ldots,a_m\}$ and
\end{itemize}
$Q=(0)$ so that $V$ is smooth and equal to $\C^m$ or $(\C^m,0)$, we
have a generic Bernstein-Sato polynomial. It does not seem
straigthforward to adapt her proof to the case where $\QQ\ne (0)$
(i.e. when $V$ is singular). Let us also say that she did not mention
the fact that the polynomial she constructed is rational even though a
detailed study of her proof shows that it is. As it appears, our main
result covers and generalizes the previous existing results in this
affine situation.

\begin{nota*}
Let $\PP$ be a prime ideal of
$\CC$. For $c$ in $\CC$, denote by $[c]_\PP$ the class of $c$ in the
quotient $\CC/\PP$ and $(c)_\PP=\frac{[c]_\PP}{1}$ this class viewed
in the fraction field of $\CC/\PP$. We naturally extend these
notations to $\CC[x]$, $\An(\CC)$ and $\CC[x,\frac{1}{f_1 \cdots f_p},
s]$.
\end{nota*}

\begin{rem*}
Using these notations, we can see that the polynomial $b(s)$ of
theorem \ref{theo:BSgen} is a Bernstein-Sato polynomial of $(f)_\PP$
for any $\PP \in V(\QQ) \smallsetminus V(h)$. This justifies the
name of a generic Bernstein-Sato polynomial on $V(\QQ)$.
\end{rem*}
As an application of theorem \ref{theo:BSgen}, we obtain some
consequences~:

\begin{cor} \label{cor1}
Fix a positive integer $d$ and a field $\k$.\\
For each $j=1,\ldots,p$,
take $f_j=\sum_{|\alpha|\le d} a_{\alpha, j} x^\alpha$ with $\alpha
\in \N^n$ and $a_{\alpha,j}$ an indeterminate. Take $a=(a_{\alpha,j})$
for $|\alpha|\le d$ and $j=1,\ldots, p$ such that we see
$f=(f_1,\ldots, f_p)$ in $\k[a][x]^p$.
Denote by $m$ the number of the $a_{\alpha,j}$'s.\\
Then there exists a finite partition of $\k^m=\cup W$ where each $W$
is a locally closed subset of $\k^m$ (i.e. $W$ is a difference of two
Zariski closed sets) such that for any $W$, there exists a polynomial
$b_W(s) \in \Q[s_1,\ldots,s_p] \smallsetminus 0$ such that for each
$a_0$ in $W$, $b_W(s)$ is in $\B^v(f(a_0,x))$.
\end{cor}

\begin{rem*}
\gauche
\begin{itemize}
\item
This corollary generalizes to the case $p\ge 2$ the main result of
A. Leykin \cite{leykin} and J. Brian\c{c}on and Ph. Maisonobe
\cite{briancon-maisonobe} in the case $p=1$.
\item
There is another way to generalize these results: Given a well
ordering $<$ on $\N^p$ compatible with sums, it is possible to prove
the existence of a partition $\k^m=\cup W$ into locally closed subsets
with the following property:
For any $W$, there exists a finite subset $G_W \subset
\k[a][x]$ such that for any $a_0 \in W$, the set $G_W(a_0)$ is a
$<$-Gr\"obner basis of the Bernstein-Sato ideal $\B^v(f(a_0,x))$, see
\cite{briancon-maisonobe} and \cite{bahloul}.
\end{itemize}
\end{rem*}

\begin{proof}[Proof of Corollary \ref{cor1}]
We remark that we can give the same statement as in corollary
\ref{cor1} for any algebraic subset $Y\subset \k^m$ as a space of
parameters. The statement of corollary \ref{cor1} will then follow
from the proof of this more general statement, that we shall give by
an induction on the dimension of $Y$. If $\dim Y=0$, the result is
trivial. Suppose $\dim Y\ge 1$. Write $Y=V_m(Q_1) \cup \cdots \cup
V_m(Q_r)$ where the $Q_i$'s are prime ideals in $\k^m$ (we identify
the maximal ideals of $\k[a]$ and the points of $\k^m$). For each $i$,
let $h_i \in \k[a] \smallsetminus Q_i$ and $b_i(s) \in \Q[s]
\smallsetminus 0$ be the $h$ and $b(s)$ of theorem \ref{theo:BSgen}
applied to $Q_i$. Now, write
\[Y=\Big( \bigcup_{i=1}^r V_m(Q_i) \smallsetminus V_m(h_i) \Big)
\bigcup Y',\]
with $Y'=\bigcup \big( V_m(Q_i) \cap V_m(h_i) \big)$ for which $\dim
Y' < \dim Y$. Apply the induction hypothesis to $Y'$. We obtain that
$Y$ is a union (not necessarily disjoint) of locally closed subsets
$V$ such that for each $V$ there exists $b_V(s) \in \Q[s]
\smallsetminus 0$ which is in $\B^v(f(a_0,x))$ for any $a_0 \in
V$. Let us show now how to obtain the annouced partition. Let $B$ be
the set of the obtained polynomials $b_V$'s. Set
$B=\{b_1,\ldots,b_e\}$. For any $i=1,\ldots,e$, let $E_i$ be the set
of the $V$'s for which $b_i=b_V$. Put
\begin{itemize}
\item $\dps W_1=\bigcup_{V\in E_1} V$,
\item $\dps W_2=\big(\bigcup_{V\in E_2} V\big) \smallsetminus \big(
  \bigcup_{V\in E_1} V\big)$,
\item[] $\vdots$
\item $\dps W_e=\big(\bigcup_{V\in E_e} V\big) \smallsetminus \big(
    \bigcup_{V\in E_1\cup \cdots \cup E_{e-1}} V\big)$.
\end{itemize}
Note that some of the $W_i$'s may be empty. The set $\{(b_1,W_1),
\ldots, (b_e,W_e)\}$ gives a partition $Y=\cup W_i$ in a way that $b_i
\in \B^v(f(a_0,x))$ for any $a_0 \in W_i$.
\end{proof}

\begin{cor} \label{cor2}
Take $f_1(a,x),\ldots,f_p(a,x) \in \mathcal{O}(U)[x]$ where
$\mathcal{O}(U)$ denotes the ring of holomorphic functions on a
open subset $U$ of $\C^m$.\\
Then there exists a finite partition of $U=\cup W$ where each $W$ is
an (analytic) locally closed subset of $U$ (i.e. each $W$ is a
difference of two analytic subsets of $U$) such that for any $W$,
there exists a rational non zero polynomial $b(s)$ which belongs to
$\B^v(f(a_0,x))$ for any $a_0 \in W$.
\end{cor}

\begin{rem*}
As it will appear in the proof, we have the same result if we
replace $\O(U)$ by $\C\{a_1,\ldots,a_m\}$ or $\k[[a_1,\ldots,a_m]]$
($\k$ being an arbitrary field).
\end{rem*} 

\begin{proof}
Let us write $f_j(a,x)=\sum g_{\alpha,j}(a)x^\alpha$ where
$g_{\alpha,j} \in \O(U)$. Let $m$ be the number of the
$g_{\alpha,j}$'s and let us introduce $m$ new variables
$b_{\alpha,j}$. Consider the (analytic) map $\phi:U \ni a \mapsto
(b_{\alpha,j}= g_{\alpha,j}(a))_{\alpha,j} \in \C^m$ where $\C$ is a
fixed arbitrary field. Now apply corollary \ref{cor1} to this
situation. Let $\k^m=\cup W$ be the obtained partition and for any
$W$, let $b_W\in \Q[s]$ be the polynomial given in \ref{cor1}. Now
apply $\phi^{-1}$. This gives a partition $U=\cup \phi^{-1}(W)$. Since
$\phi$ is analytic, the sets $\phi^{-1}(W)$ are locally closed
analytic subsets of $U$. It is then clear that for any $W$ and $a_0\in
\phi^{-1}(W)$, we have $b_W\in \B^v(f(a_0,x))$.
\end{proof}

\section*{Proof of the main theorem}

In order to prove theorem \ref{theo:BSgen}, we shall first prove the
following.

\begin{theo}\label{theo:BSglobal}
Let $\k$ be a field and $f \in \k[x]^p$. Then $\B^v(f) \cap \Q[s]$ is
not zero.
\end{theo}

Note that in \cite{briancon}, the author proved (for $p=1$) that the
global Bernstein polynomial has rational roots for any field $\k$ of
characteristic zero.
The proof of \ref{theo:BSglobal} will use the following propositions.

\begin{prop}\label{prop:subfield}
Let $\K$ be a subfield of a field $\L$. Suppose that $f\in
\K[x]^p$. Let $b(s)\in \K[s]$ be such that $b(s) f^s \in \An(\L)[s]
f^{s+v}$. Then
\[b(s)f^s \in \An(\K)[s]f^{s+v}.\]
\end{prop}
\begin{proof}
The proof is inspired by \cite{briancon} in which the case $p=1$ is
treated.
As $\L$ is a $\K$-vector space, let us take $\{1\}\cup \{l_\gamma; \,
\gamma \in \Gamma\}$ as a basis so that $\L[x,s,\frac{1} {f_1\cdots
  f_p}]f^s$ is a free $\K[x,s, \frac{1}{f_1\cdots f_p}]$-module with
$\{f^s\} \cup \{l_\gamma f^s; \, \gamma \in \Gamma\}$ as a basis.
Now let $P$ be in $\An(\L)[s]$ such that $b(s)f^s=Pf^{s+v}$.
We decompose $P=P_0+P'$ where $P_0\in \An(\K)[s]$ and $P'$ has its
coefficients in $\dps \bigoplus_{\gamma \in \Gamma} \K \cdot
l_\gamma$. Now, we have:
\[b(s)f^s=P_0 f^{s+v} + P'f^{s+v},\]
with $b(s)f^s$ and $P_0 f^{s+v}$ in $\K[x,s,\frac{1}{f_1\cdots
  f_p}]f^s$ and $P'f^{s+v}$ in $\dps \bigoplus_{\gamma \in \Gamma}
\K[x,s, \frac{1}{f_1\cdots f_p}] l_\gamma f^s$. By identification, we
obtain:
\[b(s)f^s=P_0 f^{s+v}.\]
\end{proof}

\begin{prop}(\cite{briancon} and
  \cite{briancon-maisonobe})\label{prop:loc-glob}
Given $f \in \C[x]^p$, we have~:
\begin{enumerate}
\item The set $\{\B^v_{x_0}(f) ;\, x_0 \in \C^n\}$ is finite.
\item $\B^v(f)$ is the intersection of all the $\B^v_{x_0}(f)$ where
  $x_0 \in \C^n$.
\end{enumerate}
\end{prop}

\begin{proof}[Proof of theorem \ref{theo:BSglobal}]
We shall divide the proof into two steps:\\
(a)
First, suppose that $\k=\C$. By \cite{sabbah1}, \cite{sabbah2}
and \cite{gyoja}, as mentioned in the introduction, each
$\B^v_{x_0}(f)$ contains a non zero rational polynomial. By the
previous proposition, we can take a finite product of these
polynomials and obtain a rational polynomial in $\B^v(f)$.\\
(b)
Now suppose that $\k$ is arbitrary.
Let $c_1,\ldots, c_N$ be all the coefficients that appear in the
writing of the $f_j$'s and consider the field $\K=\Q(c_1,\ldots,c_N)$.
There exist $e_1,\ldots, e_N \in \C$ and an injective morphism of
fields $\phi:\K \to \C$ such that $\phi(c_i)=e_i$ for any $i$. We
denote by the same symbol $\phi$ the natural extension of $\phi$ from
$\K[x]$ to $\C[x]$ and from $\An(\K)[s]$ to $\An(\C)[s]$.
Now, consider in $\C[s]$ the Bernstein-Sato ideal $\B^v(\phi(f))$
(where $\phi(f)=(\phi(f_1), \ldots, \phi(f_p))$. Using the result of
case (a), there exists $b(s)\in \Q[s]\smallsetminus 0$ that
belongs to $\B^v(\phi(f))$. So we have a functional equation:
\[b(s)\phi(f)^s=P \cdot \phi(f)^{s+v},\]
where $P \in \An(\C)[s]$. By proposition \ref{prop:subfield}, we can
suppose $P \in \An(\phi(\K))[s]$. Apply $\phi^{-1}$ to this
equation. Since $b(s)\in \Q[s]$, $\phi^{-1}(b(s)) =b(s)$, thus we
obtain: \[b(s) f^s=\phi^{-1}(P) \cdot f^{s+v}.\]
In conclusion $b(s)$ is in $\B^v(f)$.
\end{proof}

Now we dispose of a sufficient material to give the

\begin{proof}[Proof of theorem \ref{theo:BSgen}]
By theorem \ref{theo:BSglobal}, there exists a non zero rational
polynomial $b(s)$ in $\B^v((f)_\QQ)$. Hence, we have the following
equation:
\[b(s) \Big(\frac{[f]_\QQ}{1} \Big)^s = \frac{[U(s)]_\QQ}{[h]_\QQ}
\cdot \Big(\frac{[f]_\QQ}{1} \Big)^{s+v},\]
where $U(s) \in \An(\CC)[s]$ and $h\in \CC \smallsetminus \QQ$. It
follows that:
\[ h\, b(s) f^s - U(s) \cdot f^{s+v} \equiv 0 \mod \QQ\]
in $\CC[x,\frac{1}{f_1,\ldots,f_p},s]f^s$. Since $f_1 \cdots f_p
\notin \QQ[x]$ and $\QQ$ is prime, we obtain:
\[h \, b(s) f^s -U(s) \cdot f^{s+v} \in \QQ[x, \frac{1}{f_1 \cdots
  f_p},s] f^s.\]
\end{proof}

This article is a more general and simplified version of some results
of my thesis \cite{bahloul}.\\
\\
{\bf Acknowledgements}\\
I am grateful to my thesis advisor Michel Granger for having intoduced
me to this subject and for valuable remarks.



\end{document}